\documentclass[11pt, twoside]{article}
\usepackage{latexsym}
\usepackage{amsmath}
\usepackage{amssymb}
\usepackage[all]{xy}
\usepackage{amsfonts}
\usepackage{verbatim}
\usepackage{amsthm}
\usepackage{mathrsfs}
\usepackage{epsfig}
\usepackage{xy}
\usepackage{array}
\usepackage{stmaryrd}
\usepackage{graphicx,color}
\usepackage{xcolor}
\usepackage{tikz}
\usetikzlibrary{arrows,calc}
\usepackage{etex}
\usepackage{mathdots}
\usepackage{float}
\usepackage{graphics}
\usepackage{pdflscape}

\usepackage{anysize,hyperref}
\input xypic
\xyoption{all}
\usepackage{enumitem}
\usepackage[perpage,symbol]{footmisc}
\usepackage{setspace}
\def\A{\mathcal{A}}

\def\P{\mathcal{P}}

\def\M{\mathcal{M}}
\def\D{\mathcal{D}}
\def\C{\mathcal{C}}

\def\dr{\ar@{->}[r]}

\newcommand{\add}{\mathsf{add}\hspace{.01in}}
\def\mod{\mathsf{mod}\hspace{.01in}}
\def\proj{\mathsf{proj}\hspace{.01in}}
\def\ker{\mathsf{Ker}\hspace{.01in}}
\def\coker{\mathsf{Coker}\hspace{.01in}}
\def\End{\mathsf{End}}

\def\Hom{\mbox{Hom}}

\usepackage{bm}
\begin{document}
\baselineskip=15pt
\title{\Large{\bf Support $\bm{\tau_n}$-tilting pairs\footnotetext{This work was supported by the National Natural Science Foundation of China (Grant Nos. 11901190 and 11671221), and by the Hunan Provincial Natural Science Foundation of China (Grant No. 2018JJ3205), and by the Scientific Research Fund of Hunan Provincial Education Department (Grant No. 19B239).}}}
\medskip
\author{Panyue Zhou and Bin Zhu}

\date{}

\maketitle
\def\blue{\color{blue}}
\def\red{\color{red}}

\newtheorem{theorem}{Theorem}[section]
\newtheorem{lemma}[theorem]{Lemma}
\newtheorem{corollary}[theorem]{Corollary}
\newtheorem{proposition}[theorem]{Proposition}
\newtheorem{conjecture}{Conjecture}
\theoremstyle{definition}
\newtheorem{definition}[theorem]{Definition}
\newtheorem{question}[theorem]{Question}
\newtheorem{remark}[theorem]{Remark}
\newtheorem{remark*}[]{Remark}
\newtheorem{example}[theorem]{Example}
\newtheorem{example*}[]{Example}
\newtheorem{condition}[theorem]{Condition}
\newtheorem{condition*}[]{Condition}
\newtheorem{construction}[theorem]{Construction}
\newtheorem{construction*}[]{Construction}

\newtheorem{assumption}[theorem]{Assumption}
\newtheorem{assumption*}[]{Assumption}

\baselineskip=17pt
\parindent=0.5cm
\vspace{-6mm}

\begin{abstract}
\baselineskip=16pt
 We introduce the higher version of the notion of Adachi-Iyama-Reiten's support $\tau$-tilting pairs, which is a generalization of maximal $\tau_n$-rigid pairs in the sense of Jacobsen-J{\o}rgensen. Let $\C$ be an $(n+2)$-angulated category with an $n$-suspension
functor $\Sigma^n$ and an Opperman-Thomas cluster
tilting object. We show that relative $n$-rigid objects in $\C$ are in bijection with
$\tau_n$-rigid pairs in the $n$-abelian category $\C/\add\Sigma^n T$, and relative maximal $n$-rigid objects  in $\C$ are in bijection with
support $\tau_n$-tilting pairs.
We also show that relative $n$-self-perpendicular objects are in bijection with maximal $\tau_n$-rigid pairs.
These results generalise the work for $\C$ being $2n$-Calabi-Yau by Jacobsen-J{\o}rgensen and the work for $n=1$ by Yang-Zhu.
\\[0.2cm]
\end{abstract}

\pagestyle{myheadings}
\markboth{\rightline {\scriptsize   Panyue Zhou and Bin Zhu}}
         {\leftline{\scriptsize Support $\tau_n$-tilting pairs}}

\section{Introduction}

Higher homological algebra emerged from the higher Auslander-Reiten theory by Iyama in \cite{I1, I2}. It has gained more and more attentions since the introduction of $(n+2)$-angulated categories by Geiss-Keller-Oppermann \cite{GKO}, and the introduction of $n$-abelian categories by Jasso \cite{Ja}. The $n$-versions of triangulated or abelian categories appear as an attempt to get a better understanding of cluster tilting subcategories: some $n$-cluster tilting subcategories in triangulated categories are $(n+2)$-angulated categories and $n$-cluster tilting subcategores in abelian categories are $n$-abelian categories.

Cluster tilting objects (or subcategories) and cluster categories provide insight into cluster algebras
and their related combinatorics. They have also been used to define a new kind
of tilting theory, known as \emph{cluster tilting theory}, which generalizes APR-tilting for
hereditary algebras. Recall that the notion of cluster tilting object and related objects from \cite{BMRRT,KR,KZ,IY}.

\begin{definition}
Let $\C$ be a triangulated category with a shift functor $\Sigma$.
\begin{itemize}
\setlength{\itemsep}{1pt}
\setlength{\parsep}{1pt}
\setlength{\parskip}{1pt}
\item[\rm (1)] An object $T\in\C$ is called \emph{rigid} if ${\rm Hom}_{\C}(T, \Sigma T)=0$.

\item[\rm (2)] An object $T\in\C$ is called \emph{maximal rigid }if it is rigid and maximal with respect to the property:
$\add T=\{M\in\C\ |\ {\rm Hom}_{\C}(T\oplus M, \Sigma(T\oplus M))=0\}.$

\item[\rm (3)] An object $T\in\C$ is called \emph{cluster tilting} if
$$\add T=\{M\in\C\ |\ {\rm Hom}_{\C}(T, \Sigma M)=0\}=\{M\in\C\ |\ {\rm Hom}_{\C}(M, \Sigma T)=0\}.$$
\end{itemize}
\end{definition}
In fact, Koenig and Zhu \cite[Lemma 3.2]{KZ} show that $T$ is cluster tilting if and only if $\Hom_{\C}(T,\Sigma T)=0$ and for any object $C\in\C$, there exists a triangle $T_0\to T_1\to C\to \Sigma T_0$ where $T_0,T_1\in\add T$.

Cluster tilting theory permitted to construct abelian categories from some
triangulated categories.
By Buan-Marsh-Reiten \cite[Theorem 2.2]{BMR} in cluster categories, by Keller-Reiten
\cite[Proposition 2.1]{KR} in the $2$-Calabi-Yau case, then by Koenig-Zhu \cite[Theorem 3.3]{KZ} and Iyama-Yoshino \cite[Corollary 6.5]{IY} in the general case, one can pass from triangulated categories to abelian categories by factoring out cluster tilting subcategories.  This permits to link cluster tilting objects in triangulated categories with tilting modules in the abelian quotient categories. In fact,
Adachi, Iyama and Reiten \cite{AIR} introduced the $\tau$-tilting theory for finite dimensional algebras. It is a generalization of classical tilting theory.
They proved that for a $2$-Calabi-Yau triangulated category $\C$ with a cluster tilting object $T$, there exists a
 bijection between the basic cluster tilting objects in $\C$ and the basic support $\tau$-tilting pairs in $\mod\End_\C(T)$.
Note that each cluster tilting object is maximal rigid in a $2$-Calabi-Yau triangulated category, but
the converse is not true in general. Zhou and Zhu \cite[Theorem 2.6]{ZhZ} proved that if $\C$ is a $2$-Calabi-Yau triangulated category with a cluster tilting object, then every maximal rigid object is cluster tilting.
The bijection above in a $2$-Calabi-Yau triangulated category was generalized to any triangulated category by Yang and Zhu \cite{YZ}. They introduced the notion of relative cluster tilting objects in a triangulated category $\C$ with a cluster tilting object, which are a generalization of cluster tilting objects. Let $\C$ be a triangulated category with a cluster tilting object $T$.
They established a one-to-one correspondence between the basic relative cluster tilting objects in $\C$ and the basic support $\tau$-tilting pairs in $\mod\End_\C(T)$. For more works on this line, please see \cite{CZZ,FGL,IJY,LX, YZZ, ZZ2}.

Now we recall some higher versions of the notions and results above.

In \cite{GKO}, Geiss, Keller and Oppermann introduced $(n+2)$-angulated categories. These are
are a higher dimensional analogue of triangulated categories,
 in the sense that triangles are replaced by $(n+2)$-angles, that is, morphism sequences of length $(n+2)$. Thus a $1$-angulated category is precisely a triangulated category. They appear for example as certain cluster tilting subcategories of triangulated categories.

The notion of cluster tilting objects can be generalised to $(n+2)$-angulated categories, it dues to Oppermann and Thomas.
\begin{definition}\cite[Definition 5.3]{OT}
Let $\C$ be an $(n+2)$-angulated category with an $n$-suspension
functor $\Sigma^n$. An object $T\in\C$ is called \emph{cluster tilting} if

(1) $\Hom_{\C}(T ,\Sigma^nT)=0$.

(2)  For any object $C\in\C$, there exists an $(n+2)$-angle
$$T_0\xrightarrow{~~}T_1\xrightarrow{~~}\cdots\xrightarrow{~~}T_{n-1}\xrightarrow{~~}T_n\xrightarrow{~~}C\xrightarrow{~~}\Sigma^n T_0$$
where $T_0, T_1,\cdots,T_{n}\in\add T$.
\end{definition}

Jacobsen and J{\o}rgensen introduced the notions of maximal $n$-rigid objects and $n$-self-perpendicular objects.

\begin{definition}\cite[Definitions 0.2 and 0.3]{JJ2}
Let $\C$ be an $(n+2)$-angulated category with an $n$-suspension
functor $\Sigma^n$.
\begin{itemize}
\setlength{\itemsep}{1pt}
\setlength{\parsep}{1pt}
\setlength{\parskip}{1pt}
\item[\rm (1)] An object $T\in\C$ is called \emph{$n$-rigid} if ${\rm Hom}_{\C}(T, \Sigma^nT)=0$.

\item[\rm (2)] An object $T\in\C$ is called \emph{maximal $n$-rigid }if it is rigid and maximal with respect to the property:
$\add T=\{M\in\C\ |\ {\rm Hom}_{\C}(T\oplus M, \Sigma^n(T\oplus M))=0\}.$

\item[\rm (3)] An object $T\in\C$ is called \emph{$n$-self-perpendicular object} if
$$\add T=\{M\in\C\ |\ {\rm Hom}_{\C}(T,\Sigma^nM)=0\}=\{M\in\C\ |\ {\rm Hom}_{\C}(M,\Sigma^nT)=0\}.$$
\end{itemize}
\end{definition}

 There are also an $n$-version of the result that the quotient categories of triangulated categories by cluster tilting objects are abelian.

\begin{theorem}\emph{\cite[Theorem 0.5]{JJ1} and \cite[Theorem 3.8]{ZZ}}\label{y1}
Let $\C$ be an $(n+2)$-angulated category with an $n$-suspension
functor $\Sigma^n$ and $T$ an Opperman-Thomas cluster
tilting object with  endomorphism algebra $\Lambda=\End_{\C}(T)$.
Consider the essential image $\D$ of the functor $$\C(T,-)\colon\C \rightarrow\mod\Lambda.$$  Then $\D$ is an $n$-cluster tilting subcategory of $\mod\Lambda$. There exists a commutative diagram, as shown below, where the vertical arrow is the quotient functor and the diagonal arrow is an equivalence of categories:
$$\begin{tikzpicture}[scale=2.5]
\node (T) at (0,1) {\(\C\)};
\node (Tadd) at (0,0) {\(\C/\add\Sigma^n T\).};
\node (D) at (1,1) {\(\D\)};
\draw[->>] (T)--node[left]{} (Tadd);
\draw[->] (T)--node[above]{\(\scriptstyle{\C(T,~-)}\)} (D);
\draw[->] (Tadd)-- node[above, sloped] {\(\sim\)} (D);
\end{tikzpicture}$$
\end{theorem}

The category $\D$ is an $n$-abelian category in the sense of Jasso by \cite[Theorem 3.16]{Ja}.
By \cite[Theorem 2.3.1]{I1}, we know that $\D$ has an $n$-Auslander-Reiten
translation $\tau_n$, which is a higher analogue of the classic Auslander-Reiten translation $\tau$.
It is natural to ask if $\D$ permits a higher analogue of the $\tau$-tilting theory of \cite{AIR} for any $(n+2)$-angulated category.

In this article, we give  an affirmative answer to this question based on the work of Jacobsen and J{\o}rgensen \cite{JJ2}.

We introduce the notion of support $\tau_n$-tilting pairs which are a generalization of maximal $\tau_n$-rigid pairs.
At the same time, we also introduce the concepts of relative $n$-rigid objects, relative maximal $n$-rigid objects, relative $n$-self-perpendicular objects which can be regarded as the generalization of
$n$-rigid objects, maximal $n$-rigid objects, $n$-self-perpendicular objects, respectively.

Our main result is the following.

\begin{theorem}{\rm (}see Theorem {\rm \ref{main}} for details{\rm )}~
Let $\C$ be an $(n+2)$-angulated category with an $n$-suspension
functor $\Sigma^n$ and an Opperman-Thomas cluster
tilting object. Then there are three bijections
$$\left\{
    \begin{array}{cc}
      \mbox{isomorphism classes of } \\
      \mbox{relative $n$-rigid objects in $\C$}
    \end{array}
  \right\}
  \rightarrow
  \left\{
    \begin{array}{cc}
      \mbox{isomorphism classes of} \\
      \mbox{$\tau_n$-rigid pairs in $\D$}
    \end{array}
  \right\},$$
$$\left\{
    \begin{array}{cc}
      \mbox{isomorphism classes of } \\
      \mbox{relative maximal $n$-rigid objects in $\C$}
    \end{array}
  \right\}
  \rightarrow
  \left\{
    \begin{array}{cc}
      \mbox{isomorphism classes of} \\
      \mbox{support $\tau_n$-tilting pairs in $\D$}
    \end{array}
  \right\},$$
$$\left\{
    \begin{array}{cc}
      \mbox{isomorphism classes of } \\
      \mbox{relative $n$-self-perpendicular objects in $\C$}
    \end{array}
  \right\}
  \rightarrow
  \left\{
    \begin{array}{cc}
      \mbox{isomorphism classes of} \\
      \mbox{maximal $\tau_n$-rigid pairs in $\D$}
    \end{array}
  \right\}.$$
\end{theorem}

When $n=1$, this theorem covers \cite[Theorem 1.2]{YZ} and \cite[Theorem 1.1]{YZZ}. When $\C$ is an $2n$-Calabi-Yau $(n+2)$-angulated category, we show the following consequence which is a completion and generalisation of a recent result of Jacobsen and J{\o}rgensen \cite{JJ2}.

\begin{corollary} Let $\C$ be an $2n$-Calabi-Yau $(n+2)$-angulated category with an $n$-suspension
functor $\Sigma^n$ and an Opperman-Thomas cluster
tilting object. Then there are three bijections
$$\left\{
    \begin{array}{cc}
      \mbox{isomorphism classes of } \\
      \mbox{$n$-rigid objects in $\C$}
    \end{array}
  \right\}
  \rightarrow
  \left\{
    \begin{array}{cc}
      \mbox{isomorphism classes of} \\
      \mbox{$\tau_n$-rigid pairs in $\D$}
    \end{array}
  \right\},$$
$$\left\{
    \begin{array}{cc}
      \mbox{isomorphism classes of } \\
      \mbox{maximal $n$-rigid objects in $\C$}
    \end{array}
  \right\}
  \rightarrow
  \left\{
    \begin{array}{cc}
      \mbox{isomorphism classes of} \\
      \mbox{support $\tau_n$-tilting pairs in $\D$}
    \end{array}
  \right\},$$
$$\left\{
    \begin{array}{cc}
      \mbox{isomorphism classes of } \\
      \mbox{$n$-self-perpendicular objects in $\C$}
    \end{array}
  \right\}
  \rightarrow
  \left\{
    \begin{array}{cc}
      \mbox{isomorphism classes of} \\
      \mbox{maximal $\tau_n$-rigid pairs in $\D$}
    \end{array}
  \right\}.$$
\end{corollary}

This corollary is a generalisation of a recent result of Jacobsen and J{\o}rgensen: See Theorem \ref{main0} below.  When $n=1$, it is just the Theorem 4.1 in \cite{AIR}.

\begin{theorem}{\rm \cite[Theorem B and Theorem 1.1]{JJ2}}\label{main0}
Let $\C$ be a $2n$-Calabi-Yau $(n+2)$-angulated category with an $n$-suspension
functor $\Sigma^n$ and an Opperman-Thomas cluster
tilting object. If each indecomposable object of $\C$ is $n$-rigid, then

{\rm 1.} an object $X$ is $n$-self-perpendicular if and only if $X$ is maximal $n$-rigid;

{\rm 2.} there exists a bijection
$$\left\{
    \begin{array}{cc}
      \mbox{isomorphism classes of } \\
      \mbox{maximal $n$-rigid objects in $\C$}
    \end{array}
  \right\}
  \rightarrow
  \left\{
    \begin{array}{cc}
      \mbox{isomorphism classes of} \\
      \mbox{maximal $\tau_n$-rigid pairs in $\D$}
    \end{array}
  \right\}.$$
\end{theorem}

This article is organized as follows. In Section 2, we review some elementary definitions
and facts about $\tau$-tilting theory and $(n+2)$-angulated categories.  In Section 3, we prove our main
result and give some applications.

\section{Preliminaries}
Let $\Lambda$ be a finite dimensional $k$-algebra and $\tau$ the Auslander-Reiten translation. We denote by $\proj\Lambda$ is the category of finitely generated projective right $\Lambda$-modules.
$|M|$ denotes the number of non-isomorphic indecomposable direct summands of $M$. Support $\tau$-tilting modules were introduced by Adachi, Iyama and Reiten \cite{AIR}, which can be regarded as a generalization of tilting modules.

\begin{definition}
Let $(M, P)$ be a pair with $M\in \mod\Lambda$ and $P\in \proj\Lambda$.
\begin{itemize}
\item[(1)] $M$ is called $\tau${\rm -rigid} if {\rm Hom}$_{\Lambda}(M,\tau M)=0$.

 \item[(2)] $(M, P)$ is called a $\tau${\rm -rigid pair} if $M$ is $\tau$-rigid and {\rm Hom}$_{\Lambda}(P, M)=0$. 
\item[(3)] $(M, P)$ is called a {\rm support $\tau$-tilting pair} if it is a $\tau$-rigid pair and $|M|+|P|=|\Lambda|$.
In this case, $M$ is called a {\rm support $\tau$-tilting module}.
\end{itemize}

Adachi, Iyama and Reiten gave some equivalent characterizations of support $\tau$-tilting pairs.

\begin{remark}\cite[Corollary 2.13]{AIR}\label{rem1}
Let $(M, P)$ be a $\tau$-rigid pair. Then the following statements are equivalent:
\begin{itemize}
\item[(a)] $(M,P)$ is a support $\tau$-tilting pair.

\item[(b)] If $(M\oplus N, P)$ is a $\tau$-rigid pair for any $N\in\mod\Lambda$, then $N\in\add M$.

\item[(c)] If $\Hom_{\Lambda}(M,\tau N)=0,\Hom_{\Lambda}(N,\tau M)=0$ and $\Hom_{\Lambda}(P,N)=0$,
then $N\in\add M$.
\end{itemize}
\end{remark}
\end{definition}

We give the fact which will be used later.

\begin{lemma}{\rm \cite[Proposition 2.2]{JK}}\label{exact}
Let $\D$ be an $n$-cluster tilting subcategory of $\mod\Lambda$.

{\rm (a)}~ If~ $0\to L\to M_1\to M_2\to \cdots\to M_n$~ is
an exact sequence in $\mod\Lambda$ whose terms all lie in $\D$, then, for any $X\in\D$ there
exists the following exact sequence in $\mod\Lambda$
$$0\to{\rm Hom}_{\Lambda}(X,L)\to{\rm Hom}_{\Lambda}(X,M_1)\to{\rm Hom}_{\Lambda}(X,M_2)\to\cdots\to{\rm Hom}_{\Lambda}(X,M_n).$$

{\rm (b)}~ If~ $M_1\to M_2\to \cdots\to M_n\to N\to 0$~ is
an exact sequence in $\mod\Lambda$ whose terms all lie in $\D$, then, for any $Y\in\D$ there
exists the following exact sequence in $\mod\Lambda$
$$0\to{\rm Hom}_{\Lambda}(N,Y)\to{\rm Hom}_{\Lambda}(M_n,N)\to{\rm Hom}_{\Lambda}(M_{n-1},N)\to\cdots\to{\rm Hom}_{\Lambda}(M_1,N).$$
\end{lemma}

Jasso \cite{Ja} introduced $n$-abelian categories which are categories inhabited
by certain exact sequences with $n+2$ terms, called $n$-exact sequences. The case $n=1$
corresponds to the classical concepts of abelian categories. An important source of examples of $n$-abelian categories are $n$-cluster tilting subcategories of abelian categories. Hence it deals with $n$-cluster tilting subcategories of abelian categories.
For example, J{\o}rgensen introduced the notion of torsion classes in $n$-abelian categories, and
established the bijection between intermediate aisles and torsion classes in $n$-abelian categories associated to $n$-representation finite
algebras, see \cite[Theorem 8.5]{J}.

Recall that an $n$-abelian category in the sense of Jasso \cite[Definition 3.1]{Ja} $\M$ is \emph{projectively generated} if for every objects $M\in\M$ there exists a projective object $P\in \M$ and an epimorphism $f\colon P \to M$.
Building on work of Jasso \cite[Theorem 3.20]{Ja}, Kvamme \cite[Theorem 1.3]{K} proved that any projectively
generated $n$-abelian category is equivalent to a $n$-cluster tilting subcategory of an abelian category $\A$ with enough projectives.

We define a higher analogue of the $\tau$-tilting theory of \cite{AIR} for any $n$-abelian category.

\begin{definition}\label{tilting}
Let $\M$ be a projectively generated $n$-abelian category with an $n$-Auslander-Reiten translation $\tau_n$.
By the discussion above,  $\M$ is equivalent to a $n$-cluster tilting subcategory of an abelian category $\A$ with enough projectives. We denote by $\P$ the full subcategory of projective objects in $\A$.
Assume that $(M,P)$ is a pair with $M\in\M$ and $P\in\P$.
\begin{itemize}
\setlength{\itemsep}{1pt}
\setlength{\parsep}{1pt}
\setlength{\parskip}{1pt}
\item[(1)] An object $M\in\M$ is called $\tau_n$-rigid if $\Hom_{\A}(M, \tau_nM)=0$.

\item[(2)] $(M,P)$ is called a $\tau_n$-rigid pair in $\A$ if
$M$ is $\tau_n$-rigid and $\Hom_{\A}(P,M)=0$.

\item[(3)] $(M,P)$ is called a maximal $\tau_n$-rigid pair in $\A$ if
it satisfies:
\begin{itemize}
\setlength{\itemsep}{1pt}
\setlength{\parsep}{1pt}
\setlength{\parskip}{1pt}
\item[(i)] If $N\in\M$, then
$$N\in\add M\Longleftrightarrow
\left\{\begin{array}{l}
\Hom_{\A}(M,\tau_n N)=0,\\
\Hom_{\A}(N,\tau_nM)=0,\\
\Hom_{\A}(P,N)=0.
\end{array}
\right.$$

\item[(ii)] If $Q\in\P$, then
$$Q\in\add P\Longleftrightarrow\Hom_{\A}(Q,M)=0.$$
\end{itemize}

\item[\rm (4)]A $\tau_n$-rigid pair $(M,P)$ is called support $\tau_n$-tilting pair if it satisfies
\begin{itemize}
\item[(i)] If $(M\oplus N, P)$ is a $\tau_n$-rigid pair for any $N\in\M$, then $N\in\add M$.
\item[(ii)] If $Q\in\P$, then
$$Q\in\add P\Longleftrightarrow\Hom_{\A}(Q,M)=0.$$
\end{itemize}
\end{itemize}
\end{definition}

\begin{remark} It is obvious that a maximal $\tau_n$-rigid pair is a support $\tau_n$-tilting pair, and a support $\tau_n$-tilting pair is a $\tau_n$-rigid pair. It will be seen that all the reverse statements are not true in general, see explanation at the end of this article.\end{remark}

Let $k$ be a field and $n$ a positive integer.
In what follows, we assume that $\C$ is a $k$-linear Hom-finite $(n+2)$-angulated category with split
idempotents. The $n$-suspension functor of $\C$ is denoted by $\Sigma^n$.
We let $T$ be an Oppermann-Thomas cluster tilting object in $\C$ with  endomorphism algebra $\Lambda=\End_{\C}(T)$.
We denote by $\D$ the essential image of the functor $\C(T, -)\colon \C\to\mod\Lambda$, where $\mod\Lambda$ is the category
of finite dimensional right $\Lambda$-modules.  For any object
$M\in\C$, we denote by $\add M$ the full subcategory of $\C$ consisting of direct summands of
direct sum of finitely many copies of $M$.

\begin{remark}
Jacobsen-J{\o}rgensen \cite[Theorem 0.5]{JJ1} and Zhou-Zhu \cite[Theorem 3.8]{ZZ} show that
$\D$ is a projectively generated $n$-abelian category with an $n$-Auslander-Reiten translation $\tau_n$. Moreover,
it is equivalent to a $n$-cluster tilting subcategory of a module category $\mod\Lambda$.
In this case, the items (1), (2) and (3) in Definition \ref{tilting} are precisely define by
Jacobsen-J{\o}rgensen \cite[Definition 0.6 and Definition 0.7]{JJ2}.
Note that when $n=1$, by \cite[Proposition 2.3]{AIR} and Remark \ref{rem1}, we know that $(M, P)$ is a maximal $\tau_1$-rigid pair if and
only if it is a support $\tau$-tilting pair, and $(M, P)$ is a support $\tau_1$-rigid pair if and
only if it is a support $\tau$-tilting pair.
Hence support $\tau_n$-tilting pairs can be viewed as higher support $\tau$-tilting pairs.
\end{remark}

\begin{remark}\label{rem2}\cite[Lemma 2.1 and Lemma 2.2]{JJ1}

(1)~ The functor $\C(T, -)$ restricts to an equivalence $\add T\to\proj\Lambda$.

(2)~ $\C(T', X)\simeq \Hom_{\Lambda}(\C(T, T'),\C(T, X))$ for any $T'\in\add T$ and $X\in\C$.
\end{remark}

\medskip

The following observation is useful in the sequel.

\begin{lemma}\label{rigid}
Let $M, N\in\D$ and
$$P_n\xrightarrow{d_n}P_{n-1}\xrightarrow{d_{n-1}}\cdots\xrightarrow{d_2}P_1\xrightarrow{d_1}P_0
\xrightarrow{d_0}M\xrightarrow{~}0$$
be a minimal projective resolution of $M$. Then
${\rm Hom}_{\Lambda}(M,\tau_n N)=0$ if and only if the map
$${\rm Hom}_{\Lambda}(P_{n-1},M)\xrightarrow{~{\rm Hom}_{\Lambda}(d_n,\,N)~}{\rm Hom}_{\Lambda}(P_{n},M)$$
is surjective.  In particular, $M$ is $\tau_n$-rigid if and only if the map
$${\rm Hom}_{\Lambda}(P_{n-1},M)\xrightarrow{~{\rm Hom}_{\Lambda}(d_n,\,M)~}{\rm Hom}_{\Lambda}(P_{n},M)$$
is surjective.
\end{lemma}

\proof Since
$P_n\xrightarrow{d_n}P_{n-1}\xrightarrow{d_{n-1}}\cdots\xrightarrow{d_2}P_1\xrightarrow{d_1}P_0
\xrightarrow{d_0}M\xrightarrow{~}0$
is a minimal projective resolution of $M$, by definition, there exists a complex
$$0\xrightarrow{~}\tau_nM\xrightarrow{~}D(P_n,\Lambda)\xrightarrow{~}D(P_{n-1},\Lambda)
\xrightarrow{~}\cdots\xrightarrow{~}D(P_1,\Lambda)\xrightarrow{~}D(P_0,\Lambda),$$
where we omitted $\Hom_{\Lambda}$ because of lack of space. We claim that it is exact.
Indeed, for all $i\in\{1,2,\cdots,n-1\}$ the cohomology of the above complex at
$D\Hom_{\Lambda}(P^i,\Lambda)$ is isomorphic to $D{\rm Ext}^i_{\Lambda}(M,\Lambda)$
and hence vanishes since $M$ and $\Lambda$ belong to $\D$.

By Lemma \ref{exact}, we have
a commutative diagram of exact sequences.
$$\xymatrix{0\ar[r]&(M,\tau_nN)\ar[r]&(N,D(P^n,\Lambda))\ar[r]\ar[d]^{\wr}&(N,D(P^{n-1},\Lambda))\ar[r]\ar[d]^{\wr}&\cdots\ar[r]&(N,D(P^0,\Lambda))\ar[d]^{\wr}&&\\
&&D(P_n, N)\ar[r]^{D(d_n,\, N)\;\;}&D(P^{n-1},N)\ar[r]&\cdots\ar[r]&D(P^{0},N).}$$
Thus the assertion follows.  \qed

\begin{lemma}{\rm \cite[Lemma 2.2]{JJ2}}\label{rad}
If $M\in\C$ has no non-zero direct summands in $\add\Sigma^nT$, then there exists an
$(n+2)$-angle
$$T_0\xrightarrow{d_0}T_1\xrightarrow{d_1}T_2\xrightarrow{d_2}\cdots\xrightarrow{d_{n-1}}T_n\xrightarrow{d_n}M\xrightarrow{d_{n+1}}\Sigma^n T_0$$
where $T_0, T_1,\cdots,T_{n}\in\add T$ and for each $i\in\{1,2,\cdots,n-1\}$ the morphism $d_i\colon T_i\to T_{i+1}$ is in the radical of $\C$. Moreover,
applying the functor $\C(T, -)$ gives a complex
$$\C(T,T_0)\to\C(T,T_1)\to\C(T,T_2)\to\cdots\to\C(T,T_n)\to\C(T,M)\to 0$$
which is the start of the minimal projective resolution of $\C(T, M)$.
\end{lemma}

\section{Relative maximal $n$-rigid objects and support $\tau_n$-tilting pairs}
We first introduce the notion of relative maximal $n$-rigid objects
in $\C$, which are a generalization of maximal $n$-rigid objects.
For an object $M$ in $\C$, we use $[M](X,Y)$ to denote the subgroup of $\Hom_{\C}(X,Y)$
consisting of the morphisms from $X$ to $Y$ factoring through $\add M$.

\begin{definition}
Let $\C$ be an $(n+2)$-angulated category with with an $n$-suspension
functor $\Sigma^n$ and an Opperman-Thomas cluster tilting object.
\begin{itemize}
\item[{\rm (i)}] An object $X$ in $\C$ is called \emph{relative $n$-rigid} if there exists an Opperman-Thomas cluster tilting object $T$ such that $[\Sigma^nT](X, \Sigma^nX)=0$. In this case, $X$ is also called $\Sigma^nT$-$n$-rigid.

\item[{\rm (ii)}] An object $X$ in $\C$ is called \emph{relative maximal $n$-rigid} if there exists an Opperman-Thomas cluster tilting object $T$ such that $X$ is $\Sigma^nT$-$n$-rigid and
$$[\Sigma^nT](X\oplus M, \Sigma^d(X\oplus M))=0 \text{ implies } M\in\add X. $$
In this case, $X$ is also called maximal $\Sigma^n T$-$n$-rigid.

\item[\rm (iii)] An object $X\in\C$ is called \emph{relative $n$-self-perpendicular object} if there exists an Opperman-Thomas cluster tilting object $T$ such that
$$\add X=\{M\in\C\ |\ [\Sigma^nT](X,\Sigma^nM)=0=[\Sigma^nT](M,\Sigma^nX)\}.$$ In this case, $X$ is also called $\Sigma^n T$-$n$-self-perpendicular.
\end{itemize}
\end{definition}

\begin{remark}
From the definition,  we can immediately conclude that the following implications:
$$\begin{array}{c}
    \mbox{$\Sigma^n T$-$n$-self-perpendicular objects} \\
    \Downarrow \\
    \mbox{maximal $\Sigma^n T$-$n$-rigid objects} \\
    \Downarrow \\
    \mbox{$\Sigma^n T$-$n$-rigid objects}
  \end{array}
$$
The implications cannot be reversed in general. See explanation at the end of this article.
\end{remark}

\begin{remark}
Any $n$-rigid object in $\C$ is relative $n$-rigid.
\end{remark}

The following result indicates that $\Sigma^nT$-$n$-rigid objects and $n$-rigid objects coincide in
some cases.   We say that an $(n+2)$-angulated category $\C$ is
$2n$-Calabi-Yau if $$\C(X, \Sigma^nY)\simeq D\C(Y,\Sigma^nX)$$
naturally in $X,Y\in\C$ where $D$ is the duality functor $\Hom_k(-, k)$.

\begin{proposition}\label{prop}
If $\C$ is a $2n$-Calabi-Yau and $T$ is an Opperman-Thomas cluster tilting object in $\C$, then $X$ is $\Sigma^nT$-$n$-rigid if and only if $X$ is $n$-rigid.
In particular, $X$ is maximal $\Sigma^nT$-$n$-rigid if and only if $X$ is maximal $n$-rigid.
$X$ is $\Sigma^nT$-$n$-self-perpendicular if and only if $X$ is $n$-self-perpendicular.

\end{proposition}

\proof It is obvious that any $n$-rigid object is $\Sigma^nT$-$n$-rigid.

Now we assume that $X$ is $\Sigma^nT$-$n$-rigid.
Since $T$ is an Opperman-Thomas cluster tilting object, there exists $(n+2)$-angle
$$T_0\xrightarrow{d_0}T_1\xrightarrow{d_1}T_2\xrightarrow{d_2}\cdots\xrightarrow{d_{n-1}}T_n\xrightarrow{d_n}X\xrightarrow{d_{n+1}}\Sigma^n T_0$$
where $T_0, T_1,\cdots,T_{n}\in\add T$.
Applying the functor $\Hom_{\C}(-,\Sigma^nX)$ to the above $(n+2)$-angle, we have the following exact sequence
$$\cdots\to\Hom_{\C}(\Sigma^nT_0,\Sigma^nX)\xrightarrow{(d_{n+1},\, \Sigma^nX)}\Hom_{\C}(X,\Sigma^nX)\xrightarrow{(d_n,\, \Sigma^nX)} \Hom_{\C}(T_n,\Sigma^nX)\to \cdots.
$$
Since $X$ is $\Sigma^nT$-$n$-rigid, we obtain that $\Hom_{\C}(d_n,\Sigma^nX)$ is an injective. It follows that the morphism
$$D\Hom_{\C}(d_n,\Sigma^nX)\colon D\Hom_{\C}(T_n,\Sigma^nX)\to D\Hom_{\C}(X,\Sigma^nX)
$$
is surjective.   By the $2n$-Calabi-Yau property, we have the following commutative diagram
$$\xymatrix@C=1.2cm@R=1cm{D\Hom_{\C}(T_n,\Sigma^nX)\ar[d]^{\wr}\ar[rr]^{D {\rm Hom}_{\C}(d_n,\,\Sigma^n X)}&&D\Hom_{\C}(X,\Sigma^nX)\ar[d]^{\wr}\\
\Hom_{\C}(X,\Sigma^nT_n)\ar[rr]^{\Hom_{\C}(X,\,\Sigma^nd_n)}&&\Hom_{\C}(X,\Sigma^nX).
}
$$
Thus we get that $\Hom_{\C}(X,\,\Sigma^nd_n)\colon\Hom_{\C}(X,\Sigma^nT_n)\to \Hom_{\C}(X,\Sigma^nX)$ is surjective.
Since $X$ is $\Sigma^nT$-$n$-rigid, we have that the morphism $\Hom_{\C}(X,\Sigma^nd_n)$ is zero.
Hence $\Hom_{\C}(X,\Sigma^nX)=0$.
This shows that $X$ is $n$-rigid.  \qed

\begin{lemma}\label{lem1}
If $M,N\in\C$ has no non-zero direct summands in $\add\Sigma^nT$ and
$${\rm Hom}_{\Lambda}(\C(T,M),\tau_n\,\C(T,N))=0={\rm Hom}_{\Lambda}(\C(T,N),\tau_n\,\C(T,M)),$$
then $[\Sigma^nT](M,\Sigma^nN)=0=[\Sigma^nT](N,\Sigma^nM)$.
\end{lemma}

\proof Since $M\in\C$ has no non-zero direct summands in $\add\Sigma^nT$,  by Lemma \ref{rad},
there exists an
$(n+2)$-angle
\begin{equation}\label{u1}
T_0\xrightarrow{d_0}T_1\xrightarrow{d_1}T_2\xrightarrow{d_2}\cdots\xrightarrow{d_{n-1}}T_n\xrightarrow{d_n}M\xrightarrow{d_{n+1}}\Sigma^n T_0
\end{equation}
where $T_0, T_1,\cdots,T_{n}\in\add T$ and for each $i\in\{1,2,\cdots,n-1\}$ the morphism $d_i\colon T_i\to T_{i+1}$ is in the radical of $\C$. Moreover,
applying the functor $\C(T, -)$ gives a complex
$$\C(T,T_0)\xrightarrow{\C(T,\,d_0)}\C(T,T_1)\to\C(T,T_2)\to\cdots\to\C(T,T_n)\to\C(T,M)\to 0$$
which is the start of the minimal projective resolution of $\C(T, M)$.
By Lemma \ref{rigid}, we have that $\Hom_{\Lambda}(\C(T,\,d_0),\C(T,N))$ is surjective since ${\rm Hom}_{\Lambda}(\C(T,M),\tau_n\,\C(T,N))=0$.

Applying the functor $\C(-,\Sigma^nN)$ to the $(n+2)$-angle (\ref{u1}), we get the following exact sequence:
$$\C(\Sigma^nT_1,\Sigma^nN)\xrightarrow{~\C(\Sigma^nd_0,\, \Sigma^nN)~}\C(\Sigma^nT_0,\Sigma^nN)
\xrightarrow{~~}\C(M,\Sigma^nN)\xrightarrow{~\C(d_n,\, \Sigma^nN)~}\C(T_n,\Sigma^nN).$$
Thus we have
\begin{equation}\label{u2}
\ker\hspace{.02in}\C(d_n,\Sigma^nN)\simeq\coker\hspace{.02in}\C(\Sigma^nd_0,\,\Sigma^nN).
\end{equation}

According to Remark \ref{rem2}, we have the following commutative diagram
$$\xymatrix@C=1.5cm@R=1.5cm{\C(T_1,N)\ar[r]^-\simeq\ar[d]^{\C(d_0,\, N)} &\Hom_{\Lambda}(\C(T,T_1),\,\C(T,N))\ar[d]^{\Hom_{\Lambda}(\C(T,d_0),\,\C(T,N))} \\
\C(T_0,N) \ar[r]^-\simeq&\Hom_{\Lambda}(\C(T,T_0),\,\C(T,N))}$$
It follows that $\C(d_0, N)$ is surjective and then $\C(\Sigma^nd_0, \Sigma^nN)$ is also surjective.
Combining with (\ref{u2}), we obtain $\ker\hspace{.02in}\C(d_n,\Sigma^nN)=0$.

Now we show that $[\Sigma^nT](M,\Sigma^nN)=0$. For any morphism $u\in [\Sigma^nT](M,\Sigma^nN)$,
since $T$ is $n$-rigid, we have $ud_n=0$ implies $u\in\ker\hspace{.02in}\C(d_n,\Sigma^nM)=0$.
Hence $u=0$.  This shows that $[\Sigma^nT](M,\Sigma^nN)=0$.

Similarly, we can show $[\Sigma^nT](M,\Sigma^nN)=0$.  \qed

\begin{lemma}\label{lem2}
Assume that $M,N\in\C$ has no non-zero direct summands in $\add\Sigma^nT$. If
$$[\Sigma^nT](M,\Sigma^nN)=0=[\Sigma^nT](N,\Sigma^nM)$$
then $${\rm Hom}_{\Lambda}(\C(T,M),\tau_n\,\C(T,N))=0={\rm Hom}_{\Lambda}(\C(T,N),\tau_n\,\C(T,M)).$$
\end{lemma}

\proof
Since $M\in\C$ has no non-zero direct summands in $\add\Sigma^nT$, by Lemma \ref{rad},
there exists an
$(n+2)$-angle
\begin{equation}\label{u1}
T_0\xrightarrow{d_0}T_1\xrightarrow{d_1}T_2\xrightarrow{d_2}\cdots\xrightarrow{d_{n-1}}T_n\xrightarrow{d_n}M\xrightarrow{d_{n+1}}\Sigma^n T_0
\end{equation}
where $T_0, T_1,\cdots,T_{n}\in\add T$ and for each $i\in\{1,2,\cdots,n-1\}$ the morphism $d_i\colon T_i\to T_{i+1}$ is in the radical of $\C$. Moreover,
applying the functor $\C(T, -)$ gives a complex
$$\C(T,T_0)\xrightarrow{\C(T,\,d_0)}\C(T,T_1)\to\C(T,T_2)\to\cdots\to\C(T,T_n)\to\C(T,M)\to 0$$
which is the start of the minimal projective resolution of $\C(T, M)$.
Combining this with
Lemma \ref{rigid}, we only need to show that
$$\Hom_{\Lambda}(\C(T,\,d_0),\C(T,N))\colon \Hom_{\Lambda}(\C(T,T_1),\,\C(T,N))\to \Hom_{\Lambda}(\C(T,T_0),\,\C(T,N)) $$
 is surjective

According to Remark \ref{rem2}, we have the following commutative diagram
$$\xymatrix@C=1.5cm@R=1.5cm{\C(T_1,N)\ar[r]^-\simeq\ar[d]^{\C(d_0,\, N)} &\Hom_{\Lambda}(\C(T,T_1),\,\C(T,N))\ar[d]^{\Hom_{\Lambda}(\C(T,d_0),\,\C(T,N))} \\
\C(T_0,N) \ar[r]^-\simeq&\Hom_{\Lambda}(\C(T,T_0),\,\C(T,N))}$$
Thus it suffices for us to prove that
 $\C(d_0, N)$ is surjective. Indeed, for any morphism $v\colon T_0\to N$,
we have that $\Sigma v\circ d_{n+1}\in[\Sigma^nT](M,\Sigma^nN)=0$. So there exists a morphism $w\colon T_1\to N$
such that $\Sigma v=\Sigma w\circ\Sigma d_0$ and then $v=wd_0$.
$$\xymatrix@C=1.5cm@R=1.2cm{M\ar[r]^{d_{n+1}}\ar[r]\ar[dr]_0&\Sigma^nT_0\ar[d]^{\Sigma v}
\ar[r]^{\Sigma^nd_0}&\Sigma T_1\ar@{-->}[dl]^{\Sigma w}\\
&\Sigma N&}$$
This shows that $\C(d_0, N)$ is surjective.  

Similarly, we can show ${\rm Hom}_{\Lambda}(\C(T,N),\tau_n\,\C(T,M))=0$.  \qed
\medskip

Now we state and prove our main result.

\begin{theorem}\label{main}~
Let $\C$ be an $(n+2)$-angulated category with an $n$-suspension
functor $\Sigma^n$ and $T$ an Opperman-Thomas cluster
tilting object with  endomorphism algebra $\Lambda=\End_{\C}(T)$. Let $\D$ be the essential image of the functor $\C(T,-)\colon\C \rightarrow\mod\Lambda$, which is $n$-ableian and equivalent to the quotient category $\C/\add\Sigma^n T$. Then we have

\begin{itemize}
\item[\rm (a)]  Let $M$ be an object in $\C$ satisfying that $\add M \cap\add\Sigma^nT=\{0\}$.
Then $M$ is $\Sigma^nT$-$n$-rigid in $\C$ if and only if $\C(T, M)$ is $\tau_n$-rigid in $\D$.

\item[\rm (b)]  Decompose any object $M$ in $\C$ as $M=M_0\oplus \Sigma^nT_M$ where $\Sigma^n T_M$ is a maximal direct summand of $M$ which belongs to $\add \Sigma^nT$. Then the correspondence
$$M\longmapsto(\C(T_0, M_0),~ \C(T,T_M))$$ gives a bijection
 between the set of isomorphism classes of $\Sigma^nT$-$n$-rigid objects in $\C$
 and the set of isomorphism classes of $\tau_n$-rigid pairs in $\D$.

\item[\rm (c)]  The functor $\C(T, -)$ induces a bijection between the set of isomorphism classes
of maximal $\Sigma^nT$-$n$-rigid objects in $\C$ and the set of isomorphism
classes of support $\tau_n$-tilting pairs in $\D$.

\item[\rm (d)]  The functor $\C(T, -)$ induces a bijection between the set of isomorphism classes
of $\Sigma^nT$-$n$-self-perpendicular in $\C$ and the set of isomorphism
classes of maximal $\tau_n$-rigid pairs in $\D$.

 \end{itemize}
\end{theorem}

\proof (a) This follows from Lemma \ref{lem1} and Lemma \ref{lem2}.
\medskip

(b) We decompose any object $M$ in $\C$ as $M=M_0\oplus \Sigma^nT_M$ where $\Sigma^n T_M$ is a maximal direct summand of $M$ which belongs to $\add \Sigma^nT$.
Put $$F(M):=\big(\C(T,M_0),~\C(T,T_M)\big).$$
If $M$ is $\Sigma^nT$-$n$-rigid, according to (a), we get that
$\C(T,M_0)$ is a $\tau_n$-rigid.
Since $M$ is $\Sigma^nT$-$n$-rigid, we have $[\Sigma^nT](\Sigma^nT_M,\Sigma^nM_0)=0$ and then $\C(T_M, M_0)=0$.
By Remark \ref{rem2}, we obtain that
$\Hom_{\Lambda}\big(\C(T,T_M),~\C(T,M_0)\big)=0$.
Thus $F(M)$ is a $\tau_n$-rigid pair of $\Lambda$-modules.
On the other hand, for any $\tau_n$-rigid pair $(M, P)$ of $\Lambda$-modules, then
 $M\simeq \C(T, M')$ where $M'$ has no non-zero direct summands in $\add\Sigma^nT$ and
 $P\simeq \C(T,T')$ where $T'\in\add T$.
By definition, we have $F(M'\oplus \Sigma^nT')=(M,P)$.
It suffices to show that $M'\oplus \Sigma^nT'$ is $\Sigma^nT$-$n$-rigid,
which is a consequence of (a), Remark \ref{rem2} and the fact that
$T$ is $n$-rigid.

(c) Let $M=M_0\oplus \Sigma^nT_M$ be a maximal $\Sigma^n T$-$n$-rigid object in $\C$,
 where $T_M\in\add T$ and $M_0$ has no non-zero direct summands in $\add\Sigma^nT$.
We claim that
$$F(M):=\big(\C(T,M_0),~\C(T,T_M)\big).$$
is a support $\tau_n$-tilting pair.
Since $M$ is a maximal $\Sigma^n T$-$n$-rigid object, by (b), we know that
$F(M)$ is a $\tau_n$-rigid pair. It remains to show that $F(M)$ is a support $\tau_n$-tilting pair.

(i)~ We assume that $\big(\C(T,M_0)\oplus N,~\C(T,T_M)\big)$
is a $\tau_n$-rigid pair for any $N\in\D$.
Since $N\in\D$, we have $N\simeq\C(T,N')$ where
$N'$ has no non-zero direct summands in $\add\Sigma^nT$.
That is to say, $\big(\C(T,M_0\oplus N'),~\C(T,T_M)\big)$
is a $\tau_n$-rigid pair. By (b), we know that
$M_0\oplus N'\oplus \Sigma^nT_M=M\oplus N'$ is $\Sigma^n T$-$n$-rigid.
Since $M$ is maximal $\Sigma^n T$-$n$-rigid, we have $N'\in\add M$.
Note that $N'$ has no non-zero direct summands in $\add\Sigma^nT$,
we get that $N'\in\add M_0$ which implies $N\in\add\C(T,M_0)$.

(ii) If $Q\in\proj\Lambda$, then there exists $Q\simeq\C(T, T')$ where $T'\in\add T$.
It follows that
$$Q\in\add\C(T,T_M)\Longrightarrow\Hom_{\Lambda}(Q,\C(T,M_0))=0$$
since $\Hom_{\Lambda}(\C(T,T_M),\C(T,M_0))$=0.
Conversely, if $\Hom_{\Lambda}(Q,\C(T,M_0))=0$, by Remark \ref{rem2},
we have that $\C(T',M_0)=0$. Note that $T$ is $n$-rigid,
it is straightforward to verify that
$M\oplus \Sigma^nT'=M_0\oplus \Sigma^n T_M\oplus\Sigma^nT'$
is $\Sigma^nT$-$n$-rigid. Since $M$ is maximal $\Sigma^n T$-$n$-rigid,
we have $\Sigma^nT'\in\add M$.
Note that $M_0$ has no non-zero direct summands in $\add\Sigma^nT$,
we get that $\Sigma^nT'\in\add \Sigma^n T_M$ which implies $T'\in\add T_M$.
Thus $Q\in\add\C(T,M_0)$.

This shows that $F(M)$ is support $\tau_n$-tilting pair.

On the other hand, now assume that $(M, P)$ is a support $\tau_n$-tilting pair in $\D$.
Then $M\simeq\C(T,M_0)$ where $M_0$ has no non-zero direct summands in $\add\Sigma^nT$
and $P\simeq\C(T,T')$ where $T'\in\add T$.
By (b), we have $N:=M_0\oplus \Sigma^n T'$ is $\Sigma^n T$-$n$-rigid.
We will show that $N$ is maximal. By definition, we need to show that
if $X$ is an object such that $[\Sigma^nT](N\oplus X,\Sigma^n(N\oplus X))=0$,
then $X\in\add N$. Without loss of generality, we assume that $X$ is indecomposable.

If $X\notin\add \Sigma^nT$, then $M\oplus\C(T,X)$ is a $\tau_n$-rigid by (a).
Since $[\Sigma^nT](N\oplus X,\Sigma^n(N\oplus X))=0$, we have
$\C(T',X)=0$. By Remark \ref{rem2},
we have that $\Hom_{\Lambda}(P,\C(T,X))\simeq \C(T',X)=0$.
Thus $\big(M\oplus\C(T,X),P\big)$ is a $\tau_n$-rigid pair.
By hypothesis, $(M,P)$ is a support $\tau_n$-tilting pair,
we infer that $\C(T,X)\in\add M$ and then $X\in\add M_0\subseteq\add N$.

If $X\in\add \Sigma^nT$, then $\Sigma^{-n}X\in\add T$ implies $\C(T,\Sigma^{-n}X)\in\proj\Lambda$.
Since $[\Sigma^nT](N\oplus X,\Sigma^n(N\oplus X))=0$, we have
$\C(X,\Sigma^nM_0)=0$ and then $\C(\Sigma^{-n}X,M_0)$. By Remark \ref{rem2},
we obtain
$\Hom_{\Lambda}\big(\C(T,\Sigma^{-n}X),~\C(T,M_0)\big)\simeq \C(\Sigma^{-n}X,M_0)=0$,
By hypothesis, $(M,P)$ is a support $\tau_n$-tilting pair,
we conclude that $\C(T,\Sigma^{-n}X)\in\add \C(T,T')$.
Hence $\Sigma^{-n}X\in\add T'$ implies $X\in\add \Sigma^nT'\subseteq\add N$.
This completes the proof of (c).

(d) Assume that $M=M_0\oplus \Sigma^nT_M$ is $\Sigma^nT$-$n$-self-perpendicular in $\C$,
 where $T_M\in\add T$ and $M_0$ has no non-zero direct summands in $\add\Sigma^nT$.
We claim that
$$F(M):=\big(\C(T,M_0),~\C(T,T_M)\big)$$
is a maximal $\tau_n$-rigid pair.

(i) Suppose that $N\in\D$ satisfies
$$
\left\{\begin{array}{l}
\Hom_{\Lambda}(\C(T,M_0),\tau_n N)=0,\\
\Hom_{\Lambda}(N,\tau_n\C(T,M_0))=0,\\
\Hom_{\Lambda}(\C(T,T_M),N)=0.
\end{array}
\right.$$
For any $N\in\D$, we have $N\simeq\C(T,N_0)$ where $N_0$ has no non-zero direct summands in $\add\Sigma^nT$.\\
By Lemma \ref{lem1}, we get that $[\Sigma^nT](M_0,\Sigma^nN_0)=0=[\Sigma^nT](N_0,\Sigma^nM_0)$.

Since $\Hom_{\Lambda}(\C(T,T_M),N)=0$, by Remark \ref{rem2},
we have $\Hom_{\C}(T_M,N_0)=0$ and then $[\Sigma^nT](\Sigma^nT_M,N_0)=0$.
We infer that $[\Sigma^nT](M,\Sigma^nN_0)=[\Sigma^nT](M_0\oplus\Sigma^nT_M,\Sigma^nN_0)=0$.

Note that $[\Sigma^nT](N_0,\Sigma^n(\Sigma^nT_M))=0$ since $T$ is $n$-rigid.
Thus we obtain $$[\Sigma^nT](N_0,\Sigma^nM)=[\Sigma^nT](N_0,\Sigma^n(M_0\oplus\Sigma^nT_M)=0.$$
Since $M$ is $\Sigma^nT$-$n$-self-perpendicular, we get $N_0\in\add M$ implies
 $N_0\in\add M_0$ since $N_0$ has no non-zero direct summands in $\add\Sigma^nT$.
It follows that $N\in\add \C(T,M_0)$.
The opposite direction is obvious.

(ii)  Since $M$ is $\Sigma^nT$-$n$-self-perpendicular, thus $M$ is $\Sigma^nT$-$n$-rigid.
By (b), we know that $FM$ is a $\tau_n$-rigid pair. Hence $\Hom_{\Lambda}(\C(T,T_M),\C(T,M_0))=0$.

If $Q\in\proj\Lambda$, then there exists $Q\simeq\C(T, T')$ where $T'\in\add T$.
It follows that
$$Q\in\add\C(T,T_M)\Longrightarrow\Hom_{\Lambda}(Q,\C(T,M_0))=0,$$
since $\Hom_{\Lambda}(\C(T,T_M),\C(T,M_0))=0$.

Conversely, if $\Hom_{\Lambda}(Q,\C(T,M_0))=0$, by Remark \ref{rem2},
we have that $\C(T',M_0)=0$. Note that $T$ is $n$-rigid, we have that $[\Sigma^nT](M,\Sigma^n(\Sigma^nT'))=0$ and
it is straightforward to verify that
$M\oplus \Sigma^nT'=M_0\oplus \Sigma^n T_M\oplus\Sigma^nT'$
is $\Sigma^nT$-$n$-rigid. Thus we have $[\Sigma^nT](\Sigma^nT',\Sigma^nM)=0$.
Since $M$ is $\Sigma^nT$-$n$-self-perpendicular, we obtain $\Sigma^nT'\in\add M$.
Note that $M_0$ has no non-zero direct summands in $\add\Sigma^nT$,
we get that $\Sigma^nT'\in\add \Sigma^n T_M$ which implies $T'\in\add T_M$.
Thus $Q\in\add\C(T,M_0)$. This shows that $F(M):=\big(\C(T,M_0),~\C(T,T_M)\big)$
is a maximal $\tau_n$-rigid pair.

On the other hand, now assume that $(M, P)$ is a maximal $\tau_n$-rigid pair.
Then $M\simeq\C(T,M_0)$ where $M_0$ has no non-zero direct summands in $\add\Sigma^nT$
and $P\simeq\C(T,T')$ where $T'\in\add T$.
By (b), we have $N:=M_0\oplus \Sigma^n T'$ is $\Sigma^n T$-$n$-rigid.

We claim that $N$ is $\Sigma^nT$-$n$-self-perpendicular in $\C$.
Indeed, by definition, it suffices to show that
if $X\in\C$ is an object such that
\begin{equation}\label{relative}
[\Sigma^nT](N,\Sigma^nX)=0=[\Sigma^nT](X,\Sigma^nN),
\end{equation}
then $X\in\add N$.

Now we decompose $X=X_0\oplus \Sigma^nT_X$ where $T_X\in\add T$ and $X_0$ has no non-zero direct summands in $\add\Sigma^nT$.
By the equality (\ref{relative}), we have
$$[\Sigma^nT](M_0,\Sigma^nX_0)=0=[\Sigma^nT](X_0,\Sigma^nM_0).$$
By Lemma \ref{lem2}, we get
then $${\rm Hom}_{\Lambda}(\C(T,M_0),\tau_n\,\C(T,X_0))=0={\rm Hom}_{\Lambda}(\C(T,X_0),\tau_n\,\C(T,M_0)).$$
Since $[\Sigma^nT](N,\Sigma^nX)=0$, we have $[\Sigma^nT](\Sigma^nT',\Sigma^nX_0)=0$ implies
$\C(T',X_0)=0$.
By Remark \ref{rem2}, we obtain
$${\rm Hom}_{\Lambda}(\C(T,T'),\tau_n\,\C(T,X_0))\simeq\Hom_{\C}(T',X_0)=0.$$
Since $(M, P)$ is a maximal $\tau_n$-rigid pair, we have $\C(T,X_0)\in\add\C(T,M_0)$
implies $X_0\in\add M_0$.

Since $[\Sigma^nT](X,\Sigma^nN)=0$, we have $[\Sigma^nT](\Sigma^nT_X,\Sigma^nM_0)=0$
implies $\C(T_X,M_0)=0$.
By Remark \ref{rem2}, we obtain
$${\rm Hom}_{\Lambda}(\C(T,T_X),\tau_n\,\C(T,M_0))\simeq\Hom_{\C}(T_X,M_0)=0.$$
Since $(M, P)$ is a maximal $\tau_n$-rigid pair, we have $\C(T,T_X)\in\add\C(T,T')$
implies $T_X\in\add T'$. We also have $\Sigma^nT_X\in\Sigma^n\add T'$.

Hence $X=X_0\oplus \Sigma^nT_X\in\add(M_0\oplus\Sigma^nT')=\add N$ since $X_0\in\add M_0$ and $\Sigma^nT_X\in\Sigma^n\add T'$.
This completes the proof of (d).  \qed

\medskip

This theorem immediately yields the following  conclusion.

\begin{corollary}\label{cor1}
Let $\C$ be a $2n$-Calabi-Yau $(n+2)$-angulated category with an Opperman-Thomas cluster
tilting object.
Then there exists a bijection between the set of $n$-rigid
objects of $\C$ and the set of $\tau_n$-rigid pairs, which induces two
one-to-one correspondences to the set of maximal $n$-rigid objects and the
set of support $\tau_n$-tilting pairs,
the set of $n$-self-perpendicular objects and the
set of maximal $\tau_n$-rigid pairs.
\end{corollary}

\proof This follows from Theorem \ref{main} and Proposition \ref{prop}. \qed
\medskip

Zhou and Zhu \cite[Theorem 2.6]{ZhZ} proved that: if $\C$ is a $2$-Calabi-Yau triangulated category with a cluster tilting object, then every maximal rigid object is cluster tilting. Thus we get the following.

\begin{remark}
In Corollary \ref{cor1}, if $n=1$, it is the just the Theorem 4.1 in \cite{AIR}.
\end{remark}

\begin{remark}
In Theorem \ref{main}, if $n=1$, it covers \cite[Theorem 1.2]{YZ} and \cite[Theorem 1.1]{YZZ}.
\end{remark}
\medskip

The following result shows that maximal $\tau_n$-rigid pairs and support $\tau_n$-tilting pairs are the same under suitable conditions.

\begin{lemma}\label{lem3}
If each indecomposable object of $\C$ is $n$-rigid, then
any maximal $\tau_n$-rigid pair is precisely support $\tau_n$-tilting pair in $\D$.
\end{lemma}

\proof By definition, we know that any maximal $\tau_n$-rigid pair is a support $\tau_n$-tilting pair.
Conversely, we assume that $(M,P)$ is a support $\tau_n$-tilting pair.
If $N\in\D$ satisfies
$$
\left\{\begin{array}{l}
\Hom_{\Lambda}(M,\tau_n N)=0,\\
\Hom_{\Lambda}(N,\tau_nM)=0,\\
\Hom_{\Lambda}(P,N)=0.
\end{array}
\right.$$
Without loss of generality, we assume that $N$ is indecomposable.
Thus $N\simeq\C(T,N')$ where $N'$ has no non-zero direct summands in $\add\Sigma^nT$
and $N'$ is indecomposable. By hypothesis, $N'$ is $n$-rigid.
By Lemma \ref{lem2}, $N$ is $\tau_n$-rigid.
Hence $(M\oplus N,P)$ is a $\tau_n$-rigid pair.
Since $(M,P)$ is a support $\tau_n$-tilting pair, we get that $N\in\add M$.
The remaining conditions are clearly satisfied.
This shows that $(M,P)$ is a maximal $\tau_n$-rigid pair.  \qed
\medskip

As an application of Theorem \ref{main}, we have the following.

\begin{corollary}
In Corollary {\rm \ref{cor1}}, if each indecomposable object in $\C$ is $n$-rigid,
then Corollary {\rm \ref{cor1}} is just the Theorem B in {\rm \cite{JJ2}}.
\end{corollary}

\proof This follows from Corollary \ref{cor1}, Proposition \ref{prop} and Lemma \ref{lem3}.  \qed
\medskip

We use the following diagram to explain our main result.

Let $\C$ be an  $(n+2)$-angulated category. Jacobsen and J{\o}rgensen \cite[Theorem 1.1]{JJ2} gave the following
implications.
$$\begin{array}{c}
    \mbox{$n$-self-perpendicular objects in $\C$} \\
    \Downarrow \\
    \mbox{maximal $n$-rigid objects in $\C$} \\
    \Downarrow \\
    \mbox{$n$-rigid objects in $\C$}
  \end{array}
$$
But the implications cannot be reversed in general, see \cite[Remark 1.2]{JJ2} and a concrete example in \cite[Section 4]{JJ2}.

If $\C$ is an  $(n+2)$-angulated category with an Opperman-Thomas cluster
tilting object. We give the following implications.
$$\begin{array}{c}
    \mbox{relative $n$-self-perpendicular objects in $\C$} \\
    \Downarrow \\
    \mbox{relative maximal $n$-rigid objects in $\C$} \\
    \Downarrow \\
    \mbox{relative $n$-rigid objects in $\C$}
  \end{array}
$$
Moveover, we get the following three bijections.
$$\begin{array}{ccc}
    \mbox{relative $n$-self-perpendicular objects in $\C$} &\longleftrightarrow& \mbox{maximal $\tau_n$-rigid pairs} \\
    \Downarrow &&\Downarrow \\
    \mbox{relative maximal $n$-rigid objects in $\C$} &\longleftrightarrow& \mbox{support $\tau_n$-tilting pairs}  \\
   \Downarrow &&\Downarrow \\
    \mbox{relative $n$-rigid objects in $\C$} &\longleftrightarrow&\mbox{$\tau_n$-rigid pairs}
  \end{array}
$$

When $n=1$, the first two bijections are the same, then the bijections are the corresponding ones in \cite{YZ, YZZ}.
When $\C$ is $2n$-Calabi-Yau $(n+2)$-angulated category with an Opperman-Thomas cluster
tilting object, we obtain the following three bijections.
$$\begin{array}{ccc}
    \mbox{$n$-self-perpendicular objects in $\C$} &\longleftrightarrow& \mbox{maximal $\tau_n$-rigid pairs} \\
    \Downarrow &&\Downarrow \\
    \mbox{maximal $n$-rigid objects in $\C$} &\longleftrightarrow& \mbox{support $\tau_n$-tilting pairs}  \\
   \Downarrow &&\Downarrow \\
    \mbox{$n$-rigid objects in $\C$} &\longleftrightarrow&\mbox{$\tau_n$-rigid pairs}
  \end{array}
$$
Hence we also know that left and right implications cannot be reversed in general.

\textbf{Panyue Zhou}\\
College of Mathematics, Hunan Institute of Science and Technology, Yueyang, Hunan, 414006, People's Republic of China.\\
E-mail: \textsf{panyuezhou@163.com}\\[0.3cm]
\textbf{Bin Zhu}\\
Department of Mathematical Sciences, Tsinghua University, Beijing, 100084, People's Republic of
China.\\
E-mail: \textsf{zhu-b@mail.tsinghua.edu.cn}

\end{document}